\documentclass[11pt,a4paper]{article}

\usepackage{amsmath,amssymb}
\newtheorem{theorem}{Theorem}[section]

\newtheorem{definition}{Definition}[section]
\newtheorem{remark}{Remark}[section]
\numberwithin{equation}{section}
\newtheorem{example}{Example}[section]

\numberwithin{equation}{section}

\begin{document}

\title{\textsc{Exponential Dichotomy and Trichotomy for Skew-Evolution
Semiflows in Banach Spaces }}
\author{\textsc{Codru\c{t}a Stoica} \\
Institut de Math\' ematiques  \\
Universit\' e Bordeaux 1  \\
France  \\
e-mail: \texttt{codruta.stoica@math.u-bordeaux1.fr}
}
\date{}
\maketitle

{\footnotesize \noindent \textbf{Abstract.} The paper emphasizes the properties of exponential dichotomy and exponential trichotomy
for skew-evolution semiflows in Banach spaces, by means of
evolution semiflows and evolution cocycles. The approach is from
uniform point of view. Some characterizations which generalize
classic results are also provided.}

{\footnotesize \vspace{3mm} }

{\footnotesize \noindent \textit{Mathematics Subject
Classification:} 34D09}

{\footnotesize \vspace{2mm} }

{\footnotesize \noindent \textit{Keywords:} Evolution semiflow, evolution cocycle, skew-evolution
semiflow, exponential dichotomy, exponential trichotomy}

\section{Preliminaries}

The exponential dichotomy is one of the basic concepts in the
theory of dynamical systems and plays an important role in the
study of stable and instable manifolds. Various concepts of
dichotomy were introduced and studied by S.N. Chow and H. Leiva in
\cite{ChLe_JDE}, M. Megan, A.L. Sasu and B. Sasu in
\cite{MeSaSa_MIR}, P. Preda and C. Preda in \cite{PrPr_MIR07}. A
natural generalization of the notion of dichotomy is considered
the trichotomy, which refers at a decomposition of the space at
every moment into three closed subspaces: a stable subspace, an
instable one and a center manifold. The trichotomy was introduced
by R.J. Sacker and G.R. Sell in \cite{SaSe_JDE} and the
exponential trichotomy by S. Elaydi and O. Hajek in
\cite{ElHa_JMAA}. In recent years interesting results in the
domain of trichotomy were obtained by L.H. Popescu in \cite{Po_JMAA}, B. Sasu and A.L. Sasu
in \cite{SaSa_MZ} or L. Barreira and C. Vallis in
\cite{BaVa_JDE1}. 

A new concept of trichotomy, the null uniform exponential
trichotomy for evolution operators was introduced in
\cite{MeSt_TP05}. The study has been continued by the definition
of uniform exponential trichotomy by means of three projection
families in \cite{MeSt_IEOT}. The trichotomy is studied in the
nonuniform setting for skew-evolution semiflows in
\cite{MeSt_CAMA} and \cite{MeStBu_AADE} and for discrete
time in \cite{MeSt_TP07}. 

In our paper we extend the asymptotic
properties of exponential dichotomy and trichotomy for the newly
introduced concept of skew-evolution semiflows, which can be
considered generalization for evolution operators and skew-product semiflows.

\section{Notations. Definitions. Examples}

We consider $(X,d)$ a metric space, $V$ a Banach space and
$\mathcal{B}(V)$ the space of all bounded linear operators from
$V$ into itself. We denote the sets $T =\left\{(t,t_{0})\in
\mathbb{R}^{2}, \ t\geq t_{0}\geq 0\right\}$ and $Y=X\times V$.
Let $P:Y\rightarrow Y$ be a projector given by $P(x,v)=(x,P(x)v)$,
where $P(x)$ is a projection on $Y_{x}=\{x\}\times V$, $x\in X$.

\begin{definition}\rm\label{def_sfl_ev}
A mapping $\varphi: T\times  X\rightarrow  X$ is called
\emph{evolution semiflow} on $ X$ if following relations hold:

$(s_{1})$ $\varphi(t,t,x)=x, \ \forall (t,x)\in
\mathbb{R}_{+}\times X$

$(s_{2})$ $\varphi(t,s,\varphi(s,t_{0},x))=\varphi(t,t_{0},x),
\forall t\geq s\geq t_{0}\geq 0, x\in X$.
\end{definition}

\begin{definition}\rm\label{def_aplcoc_ev}
A mapping $\Phi: T\times  X\rightarrow \mathcal{B}(V)$ is called
\emph{evolution cocycle} over an evolution semiflow $\varphi$ if:

$(c_{1})$ $\Phi(t,t,x)=I$, the identity operator on $V$, $\forall
(t,x)\in \mathbb{R}_{+}\times X$

$(c_{2})$
$\Phi(t,s,\varphi(s,t_{0},x))\Phi(s,t_{0},x)=\Phi(t,t_{0},x),\forall
t\geq s\geq t_{0}\geq 0, x\in X$.
\end{definition}

\begin{definition}\rm\label{def_coc_ev_1}
The mapping $C: T\times Y\rightarrow Y$ defined by the relation
$C(t,s,x,v)=(\varphi(t,s,x),\Phi(t,s,x)v)$, where $\Phi$ is an
evolution cocycle over an evolution semiflow $\varphi$, is called
\emph{skew-evolution semiflow} on $Y$.
\end{definition}

\begin{example}\rm\label{ex_ses}
Let $f:\mathbb{R}_{+}\rightarrow\mathbb{R}_{+}$ be a decreasing
function on $[0,\infty)$ such that there exists
$\underset{t\rightarrow\infty}{\lim }f(t)=l$. 

Let $X$ be the closure in $\mathcal{C}(\mathbb{R}_{+},\mathbb{R})$ of the set
\[
{\{f_{t}, \ t\in \mathbb{R}_{+}\}}, \ \textrm{where} \
f_{t}(\tau)=f(t+\tau), \ \forall \tau\in \mathbb{R}_{+}.
\] 
Then the mapping 
\[
\varphi: T\times X\rightarrow  X, \
\varphi(t,s,x)=x_{t-s}
\]
is an evolution semiflow on $X$. 

Let us
consider the Banach space $V=\mathbb{R}^{p}$, $p\geq 1$, with the
norm $\left\Vert
(v_{1},...,v_{p})\right\Vert=|v_{1}|+...+|v_{p}|$. The mapping
\[
\Phi: T\times  X\rightarrow \mathcal{B}(V), \ 
\Phi(t,s,x)(v_{1},...,v_{p})=\left(
e^{\int_{s}^{t}x(\tau-s)d\tau}v_{1},...,e^{\int_{s}^{t}x(\tau-s)d\tau}v_{p}\right)
\]
is an evolution cocycle over $\varphi$ and $C=(\varphi,\Phi)$ is a
skew-evolution semiflow on $Y$.
\end{example}

\begin{definition}\rm\label{comp_2}
Two projector families $\{P_{k}\}_{k\in \{1,2\}}$ are said to be
\emph{compatible} with a skew-evolution semiflow
$C=(\varphi,\Phi)$ if

$(dc_{1})$ $P_{1}(x)+P_{2}(x)=I$,
$P_{1}(x)P_{2}(x)=P_{2}(x)P_{1}(x)=0$

$(dc_{2})$
$P_{k}(\varphi(t,s,x))\Phi(t,s,x)v=\Phi(t,s,x)P_{k}(x)v$, $k\in
\{1,2\}$

\noindent for all $t\geq s\geq t_{0}\geq 0$ and all $(x,v)\in Y.$
\end{definition}

\begin{definition}\rm\label{ued}
A skew-evolution semiflow $C=(\varphi,\Phi)$ is \emph{uniformly
exponentially dichotomic} if there exist two projector families
$\{P_{k}\}_{k\in \{1,2\}}$ compatible with $C$ and some constants
$N_{1}$, $N_{2}\geq 1$, $\nu _{1}$, $\nu _{2}>0$ such that

$(ued_{1})$ $e^{\nu _{1}(t-s)}\left\Vert
\Phi(t,t_{0},x)P_{1}(x)v\right\Vert \leq N_{1}\left\Vert
\Phi(s,t_{0},x)P_{1}(x)v\right\Vert $

$(ued_{2})$ $e^{\nu _{2}(t-s)}\left\Vert
\Phi(s,t_{0},x)P_{2}(x)v\right\Vert \leq N_{2}\left\Vert
\Phi(t,t_{0},x)P_{2}(x)v\right\Vert $

\noindent for all $t\geq s\geq t_{0}\geq 0$ and all $(x,v)\in Y.$
\end{definition}

\begin{example}\rm\label{ex_ued}
We consider $X$ and the evolution semiflow $\varphi$ as in Example
\ref{ex_ses}. Let $V=\mathbb{R}^{2}$ with the norm
$\left\Vert(v_{1},v_{2})\right\Vert=|v_{1}|+|v_{2}|$. The mapping
\[
\Phi: T\times  X\rightarrow \mathcal{B}(V), \
\Phi(t,s,x)(v)=\left(v_{1}e^{-2\int_{s}^{t}x(\tau-s)d\tau},v_{2}e^{3\int_{s}^{t}x(\tau-s)d\tau}\right)
\]
is an evolution cocycle. 

We consider the projectors
$P_{1}(x,v)=(v_{1},0)$, $P_{2}(x,v)=(0,v_{2})$. 

Then
$C=(\varphi,\Phi)$ is a uniformly exponentially dichotomic
skew-evolution semiflow with $N_{1}=N_{2}=1$, $\nu_{1}=2$,
$\nu_{2}=3$.
\end{example}

\begin{definition}\rm\label{comp_3}
Three projector families $\{P_{k}\}_{k\in \{1,2,3\}}$ are said to
be \emph{compatible} with a skew-evolution semiflow
$C=(\varphi,\Phi)$ if

$(tc_{1})$ $P_{1}(x)+P_{2}(x)+P_{3}(x)=I$,
$P_{i}(x)P_{j}(x)=P_{j}(x)P_{i}(x)=0$, $i,j\in \{1,2,3\}$, $i\neq
j$

$(tc_{2})$
$P_{k}(\varphi(t,s,x))\Phi(t,s,x)v=\Phi(t,s,x)P_{k}(x)v$, $k\in
\{1,2,3\}$

\noindent for all $t\geq s\geq t_{0}\geq 0$ and all $(x,v)\in Y.$
\end{definition}

\begin{definition}\rm\label{uet}
A skew-evolution semiflow $C$ is \emph{uniformly exponentially
trichotomic} if there exist three projector families
$\{P_{k}\}_{k\in \{1,2,3\}}$ compatible with $C$ and some
constants $N_{1}$, $N_{2}$, $N_{3}\geq 1$, $\nu_{1}$, $\nu _{2}$,
$\nu _{3}>0$ such that

$(uet_{1})$ $e^{\nu _{1}(t-s)}\left\Vert
\Phi(t,t_{0},x)P_{1}(x)v\right\Vert \leq N_{1}\left\Vert
\Phi(s,t_{0},x)P_{1}(x)v\right\Vert $

$(uet_{2})$ $e^{\nu _{2}(t-s)}\left\Vert
\Phi(s,t_{0},x)P_{2}(x)v\right\Vert \leq N_{2}\left\Vert
\Phi(t,t_{0},x)P_{2}(x)v\right\Vert $

$(uet_{3})$ $\left\Vert \Phi(t,t_{0},x)P_{3}(x)v\right\Vert \leq
N_{3}e^{\nu _{3}(t-s)}\left\Vert
\Phi(s,t_{0},x)P_{3}(x)v\right\Vert$

\hspace{10mm} $\left\Vert \Phi(s,t_{0},x)P_{3}(x)v\right\Vert \leq
N_{3}e^{\nu _{3}(t-s)}\left\Vert
\Phi(t,t_{0},x)P_{3}(x)v\right\Vert$

\noindent for all $t\geq s\geq t_{0}\geq 0$ and all $(x,v)\in Y.$
\end{definition}

\begin{example}\rm\label{ex_uet}
Let us consider function a $f$, a metric space $X$ and an
evolution semiflow $\varphi$ as in Example \ref{ex_ses}. Let
$\mu>f(0)$. 

We consider $V=\mathbb{R}^{3}$ with the norm
$\left\Vert(v_{1},v_{2},v_{3})\right\Vert=|v_{1}|+|v_{2}|+|v_{3}|$.
The mapping 
\[
\Phi: T\times  X\rightarrow \mathcal{B}(V), \
\Phi(t,s,x)(v)=
\]
\[
=(e^{-\mu(t-t_{0})+\int_{t_{0}}^{t}x(\tau-t_{0})d\tau}v_{1},\
e^{\int_{t_{0}}^{t}x(\tau-t_{0})d\tau}v_{2},\
e^{-(t-t_{0})x(0)+\int_{t_{0}}^{t}x(\tau-t_{0})d\tau}v_{3})
\]
is an evolution cocycle. 

We consider the projectors
$P_{1}(x,v)=(v_{1},0,0)$, $P_{2}(x,v)=(0,v_{2},0)$,
$P_{3}(x,v)=(0,0,v_{3})$. 

Then $C=(\varphi,\Phi)$ is uniformly
exponentially trichotomic with $N_{1}=N_{2}=N_{3}=1$,
$\nu_{1}=\mu-x(0)$, $\nu_{2}=l$, $\nu_{3} =x(0)$.
\end{example}

\begin{remark}\rm
For $P_{3}=0$ we obtain in Definition \ref{uet} the property of
uniform exponential dichotomy.
\end{remark}

\section{Main results}

To characterize the uniform exponential dichotomy we will consider
the next theorem.

\begin{theorem}\label{th_d}
A skew-evolution semiflow $C=(\varphi,\Phi)$ is uniformly
exponentially dichotomic if and only if there exist two projector
families $\{P_{k}\}_{k\in \{1,2\}}$ compatible with $C$ and a
nondecreasing function $f:[0,\infty )\rightarrow (1,\infty )$ with
the property $\underset{t\rightarrow\infty}\lim f(t)=\infty$ such
that

$(i)$ $f(t-s)\left\Vert \Phi(t,t_{0},x)P_{1}(x)v\right\Vert \leq
\left\Vert \Phi(s,t_{0},x)P_{1}(x)v\right\Vert $

$(ii)$ $f(t-s)\left\Vert \Phi(s,t_{0},x)P_{2}(x)v\right\Vert \leq
\left\Vert \Phi(t,t_{0},x)P_{2}(x)v\right\Vert $

\noindent for all $t\geq s\geq t_{0}\geq 0$ and all $(x,v)\in Y.$
\end{theorem}

\noindent \emph{Proof}. \emph{Necessity}. It is immediate if we
consider $f(t)=N^{-1}e^{\nu t}$, $t\geq 0$, where
$N=\max\{N_{1},N_{2}\}$ and $\nu=\min\{\nu_{1},\nu_{2}\}$, the
constants $N_{1}$, $N_{2}$, $\nu_{1}$, $\nu_{2}$ being given by
Definition \ref{ued}.

\emph{Sufficiency}. We will show that $(i)$ implies $(ued_{1})$.
From the definition of function $f$ there exists $\delta>0$ such that $f(\delta )>1$. We denote 
\[
\nu
=\frac{\ln f(\delta )}{\delta }>0. 
\]
Let $(t,s)\in T$. There exist
$n\in\mathbb{N}$ and $r\in[0,\delta)$ such that $t-s=n\delta+r$.
We have
\begin{equation*}
e^{\nu (t-s)}\left\Vert \Phi(t,t_{0},x)v\right\Vert \leq f(\delta
)[f(\delta )]^{n}\left\Vert \Phi(t,t_{0},x)v\right\Vert \leq
\end{equation*}%
\begin{equation*}
\leq f(\delta )[f(\delta )]^{n-1}\left\Vert \Phi(t-\delta
,t_{0},x)v\right\Vert \leq ...\leq f(\delta )\left\Vert
\Phi(t-n\delta ,t_{0},x)v\right\Vert \leq
\end{equation*}%
\begin{equation*}
\leq f(\delta )f(r)\left\Vert \Phi(t-n\delta ,t_{0},x)v\right\Vert
\leq f(\delta )\left\Vert \Phi(s,t_{0},x)v\right\Vert.
\end{equation*}
If we denote $N=f(\delta )>1$, $(ued_{1})$ follows.

By an analogous deduction we obtain that $(ii)$ implies
$(ued_{2})$.

\vspace{3mm}

Next result represent a characterization for the property of
uniform exponential trichotomy.

\begin{theorem}\label{th_t}
Let $C=(\varphi,\Phi)$ be a skew-evolution semiflow with the
property that for all $(t_{0},x,v)\in \mathbb{R}_{+}\times Y$ the
mapping $s\mapsto \left\Vert \Phi(s,t_{0},x)v\right\Vert$ is
measurable. Then $C$ is uniformly exponentially trichotomic if and
only if there exist three projector families $\{P_{k}\}_{k\in
\{1,2,3\}}$ compatible with $C$, some constants $N$, $M\geq 1$ and
a nondecreasing function $g:[0,\infty )\rightarrow (1,\infty )$
with the property $\underset{t\rightarrow\infty}\lim g(t)=\infty$
such that

$(i)$ $\left\Vert \Phi(t,t_{0},x)P_{1}(x)v\right\Vert \leq
N\left\Vert P_{1}(x)v\right\Vert $ and

\hspace{7mm} $\int_{s}^{t}\left\Vert \Phi(\tau
,t_{0},x)P_{1}(x)v\right\Vert d\tau \leq M\left\Vert
\Phi(s,t_{0},x)P_{1}(x)v\right\Vert$

$(ii)$ $\left\Vert P_{2}(x)v\right\Vert \leq N\left\Vert
\Phi(t,t_{0},x)P_{2}(x)v\right\Vert $ and

\hspace{7mm} $\int_{s}^{t}\left\Vert \Phi(\tau
,t_{0},x)P_{2}(x)v\right\Vert d\tau \leq M\left\Vert
\Phi(t,t_{0},x)P_{2}(x)v\right\Vert$

$(iii)$ $\left\Vert \Phi(t,t_{0},x)P_{3}(x)v\right\Vert \leq
g(t-s)\left\Vert \Phi(s,t_{0},x)P_{3}(x)v\right\Vert $

\hspace{7mm} $\left\Vert \Phi(s,t_{0},x)P_{3}(x)v\right\Vert \leq
g(t-s)\left\Vert \Phi(t,t_{0},x)P_{3}(x)v\right\Vert $

\noindent for all $t\geq s\geq t_{0}\geq 0$ and all $(x,v)\in Y.$
\end{theorem}

\noindent \emph{Proof}. \emph{Necessity}. It can be easily
verified. We obtain $N=\max\{N_{1},N_{2}\}$,
$M=\max\{N_{1}\nu_{1}^{-1},N_{2}\nu_{2}^{-1}\}$ and
$g(t)=N_{3}^{-1}e^{\nu_{3} t}$, the constants $N_{1}$, $N_{2}$,
$N_{3}$, $\nu_{1}$, $\nu_{2}$, $\nu_{3}$ being given by Definition
\ref{uet}.

\emph{Sufficiency}. We will prove that the relations in $(i)$
imply $(uet_{1})$. We have
\begin{equation*}
\left\Vert \Phi(t,t_{0},x)P_{1}(x)v\right\Vert = \left\Vert
\Phi(t,s,\varphi(s,t_{0},x))\Phi(s,t_{0},x)P_{1}(x)v\right\Vert\leq
\end{equation*}
\begin{equation}\label{us}
\leq N\left\Vert \Phi(s,t_{0},x)P_{1}(x)v\right\Vert ,
\end{equation}
for all $t\geq s\geq t_{0}\geq 0$ and all $(x,v)\in Y$. We
integrate the obtained relation on $[s,t]$ and we have
\begin{equation*}
(t-s)\left\Vert \Phi(t,t_{0},x)P_{1}(x)v\right\Vert \leq
N\int_{s}^{t}\left\Vert \Phi(\tau ,t_{0},x)P_{1}(x)v\right\Vert
d\tau
\end{equation*}%
\begin{equation}\label{int}
\leq MN\left\Vert \Phi(s,t_{0},x)P_{1}(x)v\right\Vert.
\end{equation}%
We obtain, according to relations (\ref{us}) and (\ref{int})
\begin{equation*}
(t-s+1)\left\Vert \Phi(t,t_{0},x)P_{1}(x)v\right\Vert \leq
N(M+1)\left\Vert \Phi(s,t_{0},x)P_{1}(x)v\right\Vert,
\end{equation*}
for all $(t,s),(s,t_{0})\in  T$ and all $(x,v)\in Y$. If we denote
\[
f(u)=\frac{u+1}{N(M+1)}, \ u\geq 0,
\]
similarly as in the proof of Theorem
\ref{th_d}, we obtain $(uet_{1})$.

An analogous deduction can be applied to prove that $(ii)$ implies
$(uet_{2})$.

To prove that the first relation in $(iii)$ implies the first
relation in $(uet_{3})$, we consider $t\geq s\geq t_{0}\geq 0$. We denote $%
n=[t-s]$. We consider $N_{3}=g(1)>1$, $\nu_{3} =\ln N_{3}>0$ and
we obtain
\begin{equation*}
\left\Vert \Phi(t,t_{0},x)P_{3}(x)v\right\Vert \leq
N_{3}\left\Vert \Phi(t-1,t_{0},x)P_{3}(x)v\right\Vert \leq ...
\end{equation*}%
\begin{equation*}
...\leq N_{3}^{n}\left\Vert \Phi(t-n,t_{0},x)P_{3}(x)v\right\Vert
\leq N_{3}^{n+1}\left\Vert \Phi(s,t_{0},x)P_{3}(x)v\right\Vert =
\end{equation*}%
\begin{equation*}
=N_{3}e^{n\nu_{3} }\left\Vert \Phi(s,t_{0},x)P_{3}(x)v\right\Vert
\leq N_{3}e^{\nu_{3} (t-s)}\left\Vert
\Phi(s,t_{0},x)P_{3}(x)v\right\Vert
\end{equation*}%
for all $t\geq s\geq t_{0}\geq 0$ and all $(x,v)\in Y$.

Similarly is obtained the second relation in $(uet_{3})$ from the
corresponding relation in $(iii)$.

\footnotesize{

}

\begin{thebibliography}{99}

\bibitem{BaVa_JDE1}
L.~Barreira and C.~Valls. Smooth center manifolds for nonuniformly
partially hyperbolic trajectories, \emph{J. Differential
Equations}, \textbf{237} (2007), 307--342.

\bibitem{ChLe_JDE} S.N.~Chow and H.~Leiva. Existence
and roughness of the exponential dichotomy for linear skew-product
semiflows in Banach spaces, \emph{J. Differential Equations},
\textbf{120} (1995), 429--477.

\bibitem{ElHa_JMAA} S.~Elaydi and O.~Hajek. Exponential
trichotomy of differential systems, \emph{J. Math. Anal. Appl.},
\textbf{129} (1988), 362--374.

\bibitem{MeSaSa_MIR} M.~Megan, A.L.~Sasu and B.~Sasu. Asymptotic Behaviour of Evolution Families, Ed. Mirton
(2003).

\bibitem{MeSt_TP05} M.~Megan and C.~Stoica. On null
uniform exponential trichotomy of evolution operators in Hilbert
spaces, \emph{Ann. Sem. Tiberiu Popoviciu on Funct. Eqs., Approx.
and Conv.}, \textbf{3} (2005), 141--150.

\bibitem{MeSt_CAMA} M.~Megan and C.~Stoica. Nonuniform
trichotomy  for skew-evolution semiflows in Banach spaces,
\emph{An. Univ. Timi\c{s}oara Ser. Mat.-Inform.}, \textbf{XLV}
Fascicle 2 (2007), 57-66.

\bibitem{MeSt_TP07} M.~Megan and C.~Stoica. Trichotomy for discrete skew-evolution semiflows in Banach
spaces, \emph{Ann. Sem. Tiberiu Popoviciu on Funct. Eqs, Approx. and
Conv.}, \textbf{5} (2007), 79-88.

\bibitem{MeSt_IEOT} \textsc{M. Megan, C. Stoica}, \textit{On uniform
exponential trichotomy of evolution operators in Banach spaces},
\emph{Integral Equations and Operators Theory},
\textbf{60} No.4 (2008), 499--506.

\bibitem{MeStBu_AADE} \textsc{M. Megan, C. Stoica, L. Buliga}, \textit{Trichotomy
for linear skew-product semiflows}, Appl. Anal. Diff. Eqs., World
Scientific (2007) 227-236

\bibitem{Po_JMAA} L.H.~Popescu. Exponential dichotomy roughness on Banach spaces,
\emph{J. Math. Anal. Appl.}, \textbf{314} (2006), 436--454.

\bibitem{PrPr_MIR07} P.~Preda, C.~Preda,
Comport\u{a}ri asimptotice ale proceselor evolutive, Ed. Mirton
(2007).

\bibitem{SaSe_JDE} R.J.~Sacker, G.R.~Sell. Existence
of dichotomies and invariant splittings for linear differential
systems III, \emph{J. Differential Equations}, \textbf{22} (1976),
497--522.

\bibitem{SaSa_MZ} B.~Sasu and A.L.~Sasu. Exponential trichotomy
and $p$-admissibility for evolution families on the real line,
\emph{Math. Z.}, \textbf{253} (2006), 515--536.

\end{thebibliography}
\end{document}